\documentclass[a4paper]{llncs}

\usepackage[latin1]{inputenc}
\usepackage[T1]{fontenc}
\usepackage{ae,aecompl}
\usepackage{amsmath,amsfonts,amssymb}
\usepackage[tableposition=bottom,font=small,labelfont=bf,format=hang]{caption}

\usepackage{booktabs}

\usepackage[dvips,colorlinks=false,bookmarks=false]{hyperref}
\hypersetup{
pdfauthor = {Lars Eirik Danielsen and Matthew G. Parker},
pdftitle = {Directed Graph Representation of Half-Rate Additive Codes over GF(4)},
pdfkeywords = {additive codes, quaternary codes, classification, graphs, circulant codes, formally self-dual codes}
}
\usepackage[dvips]{graphicx}

\DeclareMathOperator{\GF}{GF}
\DeclareMathOperator{\wt}{wt}
\DeclareMathOperator{\Aut}{Aut}
\DeclareMathOperator{\tr}{Tr}

\begin{document}

\frontmatter

\pagestyle{headings}

\mainmatter

\title{Directed Graph Representation of Half-Rate Additive Codes over $\GF(4)$}
\author{Lars Eirik Danielsen \and Matthew G. Parker}
\institute{Department of Informatics, University of Bergen, PB 7803, N-5020 Bergen, Norway\\
\texttt{\{\href{mailto:larsed@ii.uib.no}{larsed},\href{mailto:matthew@ii.uib.no}{matthew}\}@ii.uib.no}\\
\texttt{http://www.ii.uib.no/\~{}\{\href{http://www.ii.uib.no/~larsed}{larsed},\href{http://www.ii.uib.no/~matthew}{matthew}\}}
}

\maketitle

\begin{abstract}
We show that $(n,2^n)$ additive codes over $\GF(4)$ can be represented as directed graphs.
This generalizes earlier results on self-dual additive codes over $\GF(4)$, which correspond to
undirected graphs.
Graph representation reduces the complexity of code classification, and enables us to
classify additive $(n,2^n)$ codes over $\GF(4)$ of length up to 7.
From this we also derive classifications of isodual and formally self-dual codes.
We introduce new constructions of circulant and bordered circulant directed graph codes,
and show that these codes will always be isodual.
A computer search of all such codes of length up to 26 reveals that these constructions
produce many codes of high minimum distance. In particular, we find new near-extremal
formally self-dual codes of length 11 and 13, 
and isodual codes of length 24, 25, and 26 with
better minimum distance than the best known self-dual codes.
\\[\baselineskip]
\end{abstract}

\section{Introduction}

An \emph{additive} code, $\mathcal{C}$, over $\GF(4)$ of \emph{length}~$n$ is an
additive subgroup of $\GF(4)^n$. 
We denote $\GF(4) = \{0,1,\omega,\omega^2\}$, where $\omega^2 = \omega + 1$.
$\mathcal{C}$~contains $2^k$ codewords for some $0 \le k \le 2n$,
and can be defined by a $k \times n$ \emph{generator matrix}, with entries from $\GF(4)$,
whose rows span $\mathcal{C}$ additively. $\mathcal{C}$~is called an $(n,2^k)$ code.
In this paper we will only consider $(n,2^n)$, or \emph{half-rate}, codes.

The \emph{Hamming weight} of $\boldsymbol{u} \in \GF(4)^n$, denoted $\wt(\boldsymbol{u})$,
is the number of nonzero components of $\boldsymbol{u}$.
The \emph{Hamming distance} between $\boldsymbol{u}$ and $\boldsymbol{v}$
is $\wt(\boldsymbol{u} - \boldsymbol{v})$.
The \emph{minimum distance} of the code $\mathcal{C}$ is the minimal Hamming distance
between any two distinct codewords of $\mathcal{C}$. Since $\mathcal{C}$ is an additive code,
the minimum distance is also given by the smallest nonzero weight of any codeword in $\mathcal{C}$.
A code with minimum distance~$d$ is called an $(n,2^k,d)$ code.
The \emph{weight distribution} of the code $\mathcal{C}$ is the sequence
$(A_0, A_1, \ldots, A_n)$, where $A_i$ is the number of codewords of weight~$i$.
The \emph{weight enumerator} of $\mathcal{C}$ is the polynomial
\[
W_{\mathcal{C}}(x,y) = \sum_{i=0}^n A_i x^{n-i} y^i
\]

Two additive codes over $\GF(4)$, $\mathcal{C}$ and $\mathcal{C}'$, are
\emph{equivalent}~\cite{gaborit} if and only if the codewords of $\mathcal{C}$ can be mapped onto the codewords 
of $\mathcal{C}'$ by a map that consists of a permutation of coordinates (or columns of the generator matrix),
followed by multiplication of coordinates by nonzero elements from $\GF(4)$, 
followed by possible \emph{conjugation} of coordinates.
Conjugation of $x \in \GF(4)$ is defined by $\overline{x} = x^2$.
For a code of length~$n$, there is a total of $6^n n!$ such maps.
The 6 possible transformations given by scaling and conjugation of a coordinate 
are equivalent to the 6 permutations of the elements $\{1,\omega,\omega^2\}$ in the coordinate.
A map that maps a code to itself is called an \emph{automorphism} of the code.
All automorphisms of $\mathcal{C}$ make up the \emph{automorphism group}, denoted $\Aut(\mathcal{C})$.
We can use the computational algebra system \emph{Magma}~\cite{magma} to find the automorphism group 
of a code. Since, at this time, Magma has no explicit function for calculating the automorphism group of an additive code,
we use the following method, described by Calderbank~et~al.~\cite{calderbank}.
We map the $(n,2^k)$ additive code $\mathcal{C}$ over $\GF(4)$ to the $[3n,k]$ binary linear code
$\beta(\mathcal{C})$ by applying the map $0 \mapsto 000$,  $1 \mapsto 011$, $\omega \mapsto 101$,
$\omega^2 \mapsto 110$ to each generator of $\mathcal{C}$. 
We then use Magma to find $\Aut(\beta(\mathcal{C})) \cap \Aut(\beta(\GF(4)^n))$, which will be isomorphic
to $\Aut(\mathcal{C})$.

The \emph{trace map}, $\tr : \GF(4) \to \GF(2)$, is defined by
$\tr(x) = x + \overline{x}$.
The \emph{Hermitian trace inner product} of two vectors over $\GF(4)$ of length~$n$, 
$\boldsymbol{u} = (u_1,u_2,\ldots,u_n)$ and $\boldsymbol{v} = (v_1,v_2,\ldots,v_n)$, is given by
\[
\boldsymbol{u} * \boldsymbol{v} = \tr(\boldsymbol{u} \cdot \overline{\boldsymbol{v}}) 
= \sum_{i=1}^n \tr(u_i \overline{v_i}) = \sum_{i=1}^n (u_i v_i^2 + u_i^2 v_i) \pmod{2}.
\]
We define the \emph{dual} of the code $\mathcal{C}$ with respect to
the Hermitian trace inner product, $\mathcal{C}^\perp = \{ \boldsymbol{u} \in \GF(4)^n \mid
\boldsymbol{u}*\boldsymbol{c}=0 \text{ for all } \boldsymbol{c} \in \mathcal{C} \}$.
$\mathcal{C}$ is \emph{self-dual} if $\mathcal{C} = \mathcal{C}^\perp$,
\emph{formally self-dual}~\cite{hankim} if $W_{\mathcal{C}}(x,y) = W_{\mathcal{C}^\perp}(x,y)$,
and \emph{isodual} if $\mathcal{C}$ is equivalent to $\mathcal{C}^\perp$. All self-dual codes are isodual, 
all isodual codes are formally self-dual, and all formally self-dual codes are half-rate codes.
The set of \emph{linear} half-rate codes over $\GF(4)$ is a small subset of the additive half-rate codes of even length.
Optimal linear half-rate codes over $\GF(4)$ of length up to 18 were classified by
Gulliver, Östergård, and Senkevitch~\cite{gullost}.
The set of half-rate additive codes contains all self-dual, isodual, 
and formally self-dual additive codes, as well as all half-rate linear codes.

It follows from the \emph{Singleton bound}~\cite{handbook} that any half-rate additive code over $\GF(4)$ must satisfy
\[
d \le \left\lfloor\frac{n}{2}\right\rfloor + 1.
\]
$\mathcal{C}$ is called \emph{extremal} if it attains the minimum distance $d$ given by the Singleton bound, and
\emph{near-extremal} if it has minimum distance $d-1$. 
If a code has highest possible minimum distance, but is not extremal, it is called \emph{optimal}.
Han and Kim~\cite{hankim,hankim2} showed that there are no
extremal formally self-dual codes of length $n \ge 8$, and no near-extremal formally self-dual codes
of length $n=16$, $n=18$, or $n \ge 20$.
Tighter bounds on the minimum distance of self-dual additive codes over $\GF(4)$
were given by Rains and Sloane~\cite[Theorem~33]{rains}.

One of the motivations for studying self-dual additive codes over $\GF(4)$ has been
the connection to \emph{quantum error-correcting codes}~\cite{calderbank}. Non-self-dual additive
codes cannot be applied as quantum codes in the same way, but are interesting for other reasons.
Han and Kim~\cite{hankim,hankim2} studied formally self-dual additive codes over $\GF(4)$,
and showed that some of these codes have higher minimum distance than the best self-dual
codes of the same length. Additive codes may also be better than the best linear codes
of the same length. It is known that some strong binary codes can be projected onto 
additive codes over $\GF(4)$~\cite{project}. A connection between formally self-dual codes
over $\GF(4)$ and lattices has also been shown~\cite{lattices}.
We have previously studied the connection between self-dual additive codes over $GF(q^2)$,
for any prime power $q$, and \emph{weighted graphs}~\cite{weighted}. Such generalizations
could also be considered for additive codes in general, but will not be discussed in this paper.

Let $t_n$ be the number of inequivalent codes of length $n$. To find one code from each 
of the $t_n$ equivalence classes, i.e., to classify the codes of length $n$, is a hard problem.
We have previously classified all self-dual additive codes over $\GF(4)$ of length up 
to~12~\cite{selfdualgf4}, by using the fact that all such codes can be represented as
undirected graphs~\cite{bouchet,schlingemann,grassl,nest}, and that an operation called
\emph{local complementation} (LC) generates orbits of graphs that correspond to equivalence classes
of codes~\cite{bouchet,nest}.

The main result of this paper is to 
show that additive $(n,2^n)$ codes over $\GF(4)$, except for some special
cases, have representations as \emph{directed graphs}.
This basically transforms the problem of classifying such quaternary codes to a binary problem,
with reduced complexity.
We show that an algorithm by Östergård~\cite{ostergard} for checking equivalence of linear codes also
works for additive codes over $\GF(4)$.
By using this algorithm, and the fact that codes correspond to directed graphs, 
we are able to classify additive $(n,2^n)$ codes over $\GF(4)$ of length up to 7.
We find that there are more than two million inequivalent codes of length 7.
We have also checked which codes are formally self-dual, isodual, or self-dual, and give
details of this enumeration.
We introduce \emph{circulant} and \emph{bordered circulant directed graph codes}, 
and a computer search of all such codes up to length 26 reveals this
subclass of additive half-rate codes to contain many codes with high minimum distance.
Due to the structure of the generator matrices, codes from these constructions will
always be isodual, and hence also formally self-dual.
We construct new near-extremal formally self-dual codes of length 11 and 13, which
were previously unknown~\cite{hankim}. 
This also answers in the affirmative the open question of the existence of an additive 
$(13,2^{13},6)$ code~\cite{bierbrauer}.
Finally, we find isodual codes of length 24, 25, and 26 with minimum 
distance~9. The best known self-dual codes of these lengths have minimum distance 8.

\section{Directed Graph Representation}

A \emph{directed graph} is a pair $G=(V,E)$ where $V$ is a set of \emph{vertices},
and $E \subseteq V \times V$ is a set of ordered pairs called \emph{edges}. A graph with $n$ vertices
can be represented by an $n \times n$ \emph{adjacency matrix} $\Gamma$, where
$\gamma_{ij} = 1$ if $(i,j) \in E$, i.e., if there is a directed 
edge from $i$ to $j$, and $\gamma_{ij} = 0$ otherwise.
We will only consider \emph{simple} graphs, where all diagonal elements of the adjacency matrix are 0.
The special case where we always have an edge $(j,i) \in E$ whenever there is an edge $(i,j) \in E$, 
i.e., the adjacency matrix is symmetric, is called an \emph{undirected graph}.
The \emph{in-neighbourhood} of $v \in V$, denoted $NI_v \subset V$, is the set of vertices $i$ such that
there is a directed edge $(i,v) \in E$. Similarly, $NO_v \subset V$ is the \emph{out-neighbourhood} of $v$, 
i.e., the set of vertices $i$ such that there is a directed edge $(v,i)$ in $E$.
$|NI_v|$ is the \emph{indegree} of $v$, and $|NO_v|$ is the \emph{outdegree} of $v$.
Two graphs $G=(V,E)$ and $G'=(V,E')$ are \emph{isomorphic} if and only if
there exists a permutation $\pi$ of $V$ such that $(u,v) \in E \iff (\pi(u), \pi(v))
\in E'$.
A directed graph is \emph{connected}, (also known as \emph{weakly connected}), if
we can reach any vertex starting from any other vertex by traversing edges in some direction,
i.e., not necessarily in the direction they point.

\begin{definition}
A \emph{directed graph code} is an additive $(n,2^n)$ code over $\GF(4)$ that has a generator matrix
of the form $C = \Gamma + \omega I$, where $\Gamma$ is the adjacency matrix
of a simple directed graph and $I$ is the identity matrix.
\end{definition}

\begin{proposition}\label{dualtranspose}
Given a directed graph code $\mathcal{C}$ with generator matrix $C = \Gamma + \omega I$, 
its dual $\mathcal{C}^\perp$ is generated by $C^T$.
\end{proposition}
\begin{proof}
We must show that for any $c \in \mathcal{C}$ and any $c' \in \mathcal{C}^\perp$, the trace inner product 
$c * c' = \tr(c \cdot \overline{c'}) = 0$. Let $c = aC$ and $c' = bC^T$, with $a,b \in \GF(2)^n$.
Then $c * c' = (aC) * (bC^T) = \tr( (aC) \cdot \overline{(bC^T)} ) = \tr( (aC)(b\overline{C}^T)^T )
= \tr( aC \overline{C}b^T )$, which must be 0 if all elements of $C\overline{C}$ are from $\GF(2)$.
This is clearly the case, since $C\overline{C} = (\Gamma + \omega I)(\Gamma + \omega^2 I) = \Gamma^2 + \Gamma + I$.\qed
\end{proof}

\begin{theorem}\label{graphthm}
Given an additive $(n,2^n)$ code $\mathcal{C}$ over $\GF(4)$ with generator matrix $C$,
there exists a directed graph code equivalent to $\mathcal{C}$, except when $C=A+\omega B$ is such
that all $2^n$ sets of $n$ columns $(\{\boldsymbol{a}_1,\boldsymbol{b}_1\}, 
\{\boldsymbol{a}_2,\boldsymbol{b}_2\}, \ldots, \{\boldsymbol{a}_n,\boldsymbol{b}_n\})$ of $(A \mid B)$,
where $\{\boldsymbol{a}_i, \boldsymbol{b}_i\}$ means that we choose either $\boldsymbol{a}_i$ 
or $\boldsymbol{b}_i$, are linearly dependent.
\end{theorem}
\begin{proof}
We can write $C=A+\omega B$, with $(A \mid B)$ a binary $n \times 2n$ matrix. From
the fact that the rows of $C$ additively span a vector space of dimension $n$, it follows that 
$(A \mid B)$ has full rank. If the $n \times n$ submatrix $B$ also has full rank over $\GF(2)$, 
we simply perform the basis change $B^{-1} (A \mid B) = (\Gamma' | I)$. Any non-zero elements on
the diagonal of $\Gamma'$ can simply be set to zero, effected by conjugating the corresponding coordinates
of $\Gamma' + \omega I$, to obtain an equivalent directed graph code generated by $\Gamma + \omega I$.

In the case that $B$ does not have full rank, we must show that there is a code $\mathcal{C'}$, 
equivalent to $\mathcal{C}$, generated by $A' + \omega B'$ where $B'$ does have full rank. 
Then we can apply the method described in the first part of this proof to obtain the graph form. 
Let the columns of $A$ be denoted $(\boldsymbol{a}_1, \boldsymbol{a}_2, \ldots, \boldsymbol{a}_n)$ 
and the columns of $B$ be denoted $(\boldsymbol{b}_1, \boldsymbol{b}_2, \ldots, \boldsymbol{b}_n)$. 
Observe that multiplying column $i$ of $C$ by $\omega^2$, followed by 
conjugation of the same column, has the effect of swapping columns $\boldsymbol{a}_i$ and $\boldsymbol{b}_i$ 
in $(A \mid B)$. Since we know that at least one of the $2^n$ possible choices 
$(\{\boldsymbol{a}_1,\boldsymbol{b}_1\}, \{\boldsymbol{a}_2,\boldsymbol{b}_2\}, \ldots, 
\{\boldsymbol{a}_n,\boldsymbol{b}_n\})$ span a vector space of dimension $n$, we can find a 
matrix $(A' \mid B')$ where $B'$ has full rank.\qed
\end{proof}

It follows from Prop.~\ref{dualtranspose} that a directed graph code is self-dual if and only if
its generator matrix is symmetric, i.e., it is in fact an undirected graph code.
The fact that all self-dual additive codes over $\GF(4)$ can be represented as 
undirected graphs is well known~\cite{bouchet,schlingemann,grassl,nest},
and was used to classify all self-dual additive codes up to length 12~\cite{selfdualgf4}.
Theorem~\ref{graphthm} is a generalization of this result to the much larger classes
of directed graphs and half-rate additive codes over $\GF(4)$.
As stated in Theorem~\ref{graphthm}, some special codes do not have graph representations.
These are codes that will typically not be of interest, such as codes with a
generator matrix that contains an all-zero column or a set of linearly dependent binary columns
(up to scaling by $\omega$ or $\omega^2$).

\begin{figure}
\centering
\includegraphics[width=.45\linewidth]{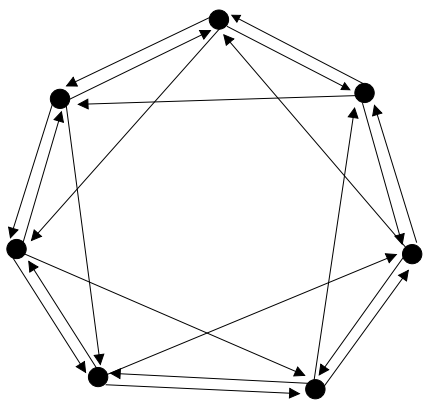}
\caption{Directed Graph Representation of a $(7,2^7,4)$ Code}\label{figstar}
\end{figure}

\begin{example}
We consider an additive $(7,2^7,4)$ code, $\mathcal{C}$,
generated by
\[
C = \left(\!\!\begin{array}{ccccccc}
\omega^2 & 0      & 0      & 1      & 0        & \omega   & \omega  \\
1        & 0      & 0      & 0      & 1        & 1        & 1       \\
0        & 0      & \omega & 0      & \omega^2 & \omega^2 & \omega  \\
0        & 0      & 1      & 1      & \omega^2 & \omega   & \omega^2 \\
1        & \omega & 0      & 1      & \omega^2 & \omega^2 & 0       \\
1        & 1      & 1      & 1      & \omega^2 & 1        & 1       \\
0        & 0      & 1      & \omega & 1        & \omega^2 & 1
\end{array}\!\!\right) = A + \omega B,
\]
\[
(A \mid B) = \left(\!\!\begin{array}{ccccccc|ccccccc}
1 & 0 & 0 & 1 & 0 & 0 & 0     & 1 & 0 & 0 & 0 & 0 & 1 & 1 \\
1 & 0 & 0 & 0 & 1 & 1 & 1     & 0 & 0 & 0 & 0 & 0 & 0 & 0 \\
0 & 0 & 0 & 0 & 1 & 1 & 0     & 0 & 0 & 1 & 0 & 1 & 1 & 1 \\
0 & 0 & 1 & 1 & 1 & 0 & 1     & 0 & 0 & 0 & 0 & 1 & 1 & 1 \\
1 & 0 & 0 & 1 & 1 & 1 & 0     & 0 & 1 & 0 & 0 & 1 & 1 & 0 \\
1 & 1 & 1 & 1 & 1 & 1 & 1     & 0 & 0 & 0 & 0 & 1 & 0 & 0 \\
0 & 0 & 1 & 0 & 1 & 1 & 1     & 0 & 0 & 0 & 1 & 0 & 1 & 0
\end{array}\!\!\right).
\]
We swap column $a_6$ with $b_6$ and column $a_7$ with $b_7$ to get the matrix
\[
(A' \mid B') = \left(\!\!\begin{array}{ccccccc|ccccccc}
1 & 0 & 0 & 1 & 0 & 1 & 1     & 1 & 0 & 0 & 0 & 0 & 0 & 0 \\
1 & 0 & 0 & 0 & 1 & 0 & 0     & 0 & 0 & 0 & 0 & 0 & 1 & 1 \\
0 & 0 & 0 & 0 & 1 & 1 & 1     & 0 & 0 & 1 & 0 & 1 & 1 & 0 \\
0 & 0 & 1 & 1 & 1 & 1 & 1     & 0 & 0 & 0 & 0 & 1 & 0 & 1 \\
1 & 0 & 0 & 1 & 1 & 1 & 0     & 0 & 1 & 0 & 0 & 1 & 1 & 0 \\
1 & 1 & 1 & 1 & 1 & 0 & 0     & 0 & 0 & 0 & 0 & 1 & 1 & 1 \\
0 & 0 & 1 & 0 & 1 & 1 & 0     & 0 & 0 & 0 & 1 & 0 & 1 & 1
\end{array}\!\!\right),
\]
where $B'$ has full rank. We can then obtain the matrix
\[
(\Gamma' \mid I) = {B'}^{-1}(A' \mid B') = \left(\!\!\begin{array}{ccccccc|ccccccc}
1 & 0 & 0 & 1 & 0 & 1 & 1     & 1 & 0 & 0 & 0 & 0 & 0 & 0 \\
0 & 0 & 1 & 0 & 1 & 0 & 1     & 0 & 1 & 0 & 0 & 0 & 0 & 0 \\
1 & 0 & 1 & 1 & 1 & 0 & 0     & 0 & 0 & 1 & 0 & 0 & 0 & 0 \\
1 & 0 & 1 & 0 & 0 & 1 & 0     & 0 & 0 & 0 & 1 & 0 & 0 & 0 \\
0 & 1 & 1 & 1 & 0 & 0 & 0     & 0 & 0 & 0 & 0 & 1 & 0 & 0 \\
1 & 1 & 0 & 0 & 0 & 1 & 1     & 0 & 0 & 0 & 0 & 0 & 1 & 0 \\
0 & 1 & 0 & 0 & 1 & 1 & 1     & 0 & 0 & 0 & 0 & 0 & 0 & 1
\end{array}\!\!\right).
\]
By setting the diagonal of $\Gamma'$ to zero, we get the adjacency matrix of a simple directed graph,
\[
\Gamma = \left(\!\!\begin{array}{ccccccc}
0 & 0 & 0 & 1 & 0 & 1 & 1 \\
0 & 0 & 1 & 0 & 1 & 0 & 1 \\
1 & 0 & 0 & 1 & 1 & 0 & 0 \\
1 & 0 & 1 & 0 & 0 & 1 & 0 \\
0 & 1 & 1 & 1 & 0 & 0 & 0 \\
1 & 1 & 0 & 0 & 0 & 0 & 1 \\
0 & 1 & 0 & 0 & 1 & 1 & 0
\end{array}\!\!\right).
\]
This graph is shown in Fig.~\ref{figstar}. 
$\Gamma + \omega I$ generates a $(7,2^7,4)$ directed graph code equivalent to $\mathcal{C}$.
\end{example}

\section{Classification}

Since we have shown in Theorem~\ref{graphthm} that, except for some special cases,
additive codes over $\GF(4)$ are equivalent to directed graph codes,
it follows that to classify such codes, we only need to consider directed graphs.
Furthermore, we only need to consider non-isomorphic graphs, since if two directed graphs
are isomorphic, the corresponding directed graph codes are equivalent by a permutation of coordinates.
All non-isomorphic connected directed graphs on up to 7 vertices can be generated in a few hours by using
tools provided with the software package \emph{nauty}~\cite{nauty}. For an enumeration of
these graphs, see sequence A003085 in \emph{The On-Line Encyclopedia of Integer Sequences}~\cite{sequences}.

Connected graphs correspond to \emph{indecomposable} codes.
A code is decomposable if it can be written as the \emph{direct sum} of two smaller codes.
For example, let $\mathcal{C}$ be an $(n,2^n,d)$ code and $\mathcal{C}'$ an $(n',2^{n'},d')$ code. The
direct sum, $\mathcal{C} \oplus \mathcal{C}' = \{u||v \mid u \in \mathcal{C}, v \in \mathcal{C}'\}$,
where $||$ means concatenation, is an $({n+n'},2^{n+n'},\min\{d,d'\})$ code.
It follows that all decomposable codes of length~$n$ can be classified easily once
all indecomposable codes of length less than $n$ are known.
The total number of codes of length $n$, $t_n$, is easily derived from the numbers $i_n$ of
indecomposable codes by using the \emph{Euler transform}~\cite{sloane2},
\begin{eqnarray*}
c_n &=& \sum_{d|n} d i_d\\
t_1 &=& c_1\\
t_n &=& \frac{1}{n}\left( c_n + \sum_{k=1}^{n-1} c_k t_{n-k} \right).
\end{eqnarray*}

\begin{figure}
\centering
\includegraphics[width=.5\linewidth]{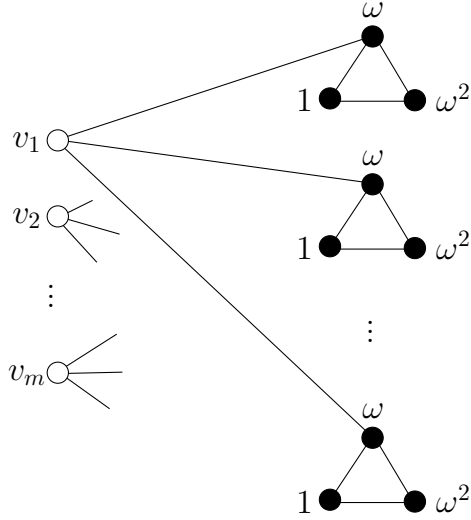}
\caption{Example of Code Graph for Checking Equivalence}\label{figcanon}
\end{figure}

To check whether two additive codes over $\GF(4)$ are equivalent, we use a modified version
of an algorithm originally devised by Östergård~\cite{ostergard} for checking equivalence of linear codes.
We show that this method also works for additive codes. 
An additive code over $\GF(4)$ is mapped to an undirected colored \emph{code graph} in 
the following way. (Note that this representation is not related to the directed graph 
representation defined previously.)
First, we find a set of vectors of some weights
that generate the code. Often, the set of all vectors of minimum weight $d$ will suffice, 
otherwise, we add all vectors of weight $d+1$, and then all vectors of weight $d+2, \ldots$, as necessary.
For each vector $c_i$ in the resulting set, add a vertex $v_i$ to the code graph. 
Also add $n$ sets of three vertices, where $n$ is the length of the code. 
Each set of three vertices represent the non-zero elements $\{1,\omega,\omega^2\} \in \GF(4)$ in one 
coordinate. 
In every set, each of the three vertices is connected to each of the two other
by undirected edges, to form a cycle. 
(This corresponds to the fact that any permutation of the symbols $\{1,\omega,\omega^2\}$
in each coordinate of the code gives an equivalent code.) Let the vertices $v_i$ have one color, 
and the $3n$ other vertices have a different color. 
Add edges between vertex $v_i$ and the $n$ 3-cycles corresponding to the codeword $c_i$.
E.g., if $c_i$ has $\omega$ in coordinate $j$, then there is an edge between $v_i$ and the element 
labelled $\omega$ in the $j$th 3-cycle. As an example, Fig.~\ref{figcanon} shows the case
where $c_1 = (\omega, \omega, \ldots, \omega)$.
The resulting code graph is then \emph{canonized}, i.e., relabelled, but with coloring preserved,
using the \emph{nauty} software~\cite{nauty}. If two graphs are isomorphic, their canonical representations are guaranteed
to be the same. Hence, if two codes are equivalent, their canonical code graphs will be identical.
Furthermore, as an alternative to the method described in the introduction, we can find
the automorphism group of a code as the automorphism group of its code graph,
i.e., the set of all vertex permutations that map the code graph to itself.

To classify codes of length up to 7, we take all non-isomorphic connected directed graphs, map
them to codes, and canonize the corresponding code graphs as described above. 
All duplicates are removed to obtain one representative from each equivalence class.
The special form of the generator matrix of a directed graph code makes it easier to find 
all codewords of small weight. If $\mathcal{C}$ is generated by $C = \Gamma + \omega I$,
then any codeword formed by adding $i$ rows of $C$ must have weight at least $i$.
This means that we can find all codewords of weight $i$ by only considering sums 
of at most $i$ rows of $C$. This property also helps when we want to find
the minimum distance of a code. Furthermore, if we wanted to exclude codes with minimum distance one from
our classification, it would suffice to exclude graphs where some vertex has outdegree zero, since this would
imply that there is a row in the generator matrix with weight one.

We also note another special property of directed graph codes:
Given a directed graph code, $\mathcal{C}$, with generator matrix $\Gamma + \omega I$,
it can be verified that the additive code over $\mathbb{Z}_4$ given by $2\Gamma + I$ always has
the same weight distribution as $\mathcal{C}$. 
We may therefore replace the elements from $\GF(4)$ with elements from 
$\mathbb{Z}_4$ by the mapping $0 \mapsto 0$, $1 \mapsto 2$, $\omega \mapsto 1$, $\omega^2 \mapsto 3$.

Table~\ref{tab:all} gives the number of half-rate additive codes over $\GF(4)$ by length and minimum distance.
Note that only indecomposable codes are counted, and that the special cases in Theorem~\ref{graphthm}, 
that do not have a directed graph representation, are not included.
A database containing one representative from each equivalence class is 
available at \url{http://www.ii.uib.no/~larsed/directed/}.
Table~\ref{tab:fsd} and Table~\ref{tab:isodual} give the numbers of formally self-dual and isodual codes.
(Note that the 240 formally self-dual $(6,2^6,3)$-codes and 3 formally self-dual $(7,2^7,4)$-codes
were also found by Han and Kim~\cite{hankim}.)
For completeness, we include in Table~\ref{tab:selfdual} the number of self-dual codes, although
we have previously classified these up to length~12~\cite{selfdualgf4}.

\begin{table}[p]
\centering
\caption{Number of Half-Rate Additive Codes over $\GF(4)$}
\label{tab:all}
\begin{tabular}{crrrrrr}
\toprule
$d \backslash n$ & 2 & 3 &  4 &   5 &      6 &          7 \\
\midrule
1                & 1 & 4 & 27 & 322 &   8509 &    686,531 \\
2                & 1 & 3 & 21 & 262 &   9653 &  1,279,641 \\
3                &   &   &  1 &   9 &    644 &    253,635 \\
4                &   &   &    &     &      1 &          3 \\
\midrule
Total            & 2 & 7 & 49 & 593 & 18,807 &  2,219,810 \\
\bottomrule
\end{tabular}
\end{table}

\begin{table}[p]
\centering
\caption{Number of Formally Self-Dual Additive Codes over $\GF(4)$}
\label{tab:fsd}
\begin{tabular}{crrrrrr}
\toprule
$d \backslash n$ & 2 & 3 &  4 &   5 &      6 &        7 \\
\midrule
1                & 1 & 1 & 10 &  55 &   1082 &   36,129 \\
2                & 1 & 2 & 12 &  79 &   2348 &  192,201 \\
3                &   &   &  1 &   5 &    240 &   55,711 \\
4                &   &   &    &     &      1 &        3 \\
\midrule
Total            & 2 & 3 & 23 & 139 &   3671 &  284,044 \\
\bottomrule
\end{tabular}
\end{table}

\begin{table}[p]
\centering
\caption{Number of Isodual Additive Codes over $\GF(4)$}
\label{tab:isodual}
\begin{tabular}{crrrrrr}
\toprule
$d \backslash n$ & 2 & 3 &  4 &   5 &      6 &       7 \\
\midrule
1                & 1 & 1 &  8 &  27 &    344 &    3243 \\
2                & 1 & 2 & 10 &  45 &    598 &    8517 \\
3                &   &   &  1 &   5 &    124 &    3299 \\
4                &   &   &    &     &      1 &       3 \\
\midrule
Total            & 2 & 3 & 19 &  77 &   1067 &  15,062 \\
\bottomrule
\end{tabular}
\end{table}

\begin{table}[p]
\centering
\caption{Number of Self-Dual Additive Codes over $\GF(4)$}
\label{tab:selfdual}
\begin{tabular}{crrrrrr}
\toprule
$d \backslash n$ & 2 & 3 &  4 &   5 &      6 &    7 \\
\midrule
2                & 1 & 1 &  2 &   3 &      9 &   22 \\
3                &   &   &    &   1 &      1 &    4 \\
4                &   &   &    &     &      1 &      \\
\midrule
Total            & 1 & 1 &  2 &   4 &     11 &   26 \\
\bottomrule
\end{tabular}
\end{table}

\section{Circulant Directed Graph Codes}

Since it is infeasible to study all half-rate additive codes of lengths much higher than those
classified in the previous section, we restrict our search space to codes corresponding to
graphs with \emph{circulant} adjacency matrices. A matrix is circulant if the $i$th row is equal to the
first row, cyclically shifted $i-1$ times to the right. 
The generator matrix of a \emph{directed graph code} is obtained by setting all diagonal elements of the
circulant adjacency matrix to $\omega$. There are $2^{n-1}$ such codes of length $n$, some of which may be equivalent.
We also consider \emph{bordered circulant} adjacency matrices: Given a length $n$ circulant graph code with
generator matrix $C$, we obtain a code of length $n+1$ with generator matrix
\[
\begin{pmatrix}
\omega & 1 & \cdots & 1 \\
1      &   &        &   \\
\vdots &   &   C    &   \\
1      &   &        &   \\
\end{pmatrix}.
\]
There are $2^{n-2}$ such codes of length $n$, some of which may be equivalent.
For each $n$ up to 26, we have counted, up to equivalence, all circulant and bordered circulant directed 
graph codes.
The number of codes of the highest found minimum distance for each $n$ is given in Table~\ref{tab:circulant}.
A database of these codes is available at \url{http://www.ii.uib.no/~larsed/directed/}.

\begin{table}
\centering
\caption{Number of Circulant and Bordered Circulant Directed Graph Codes of Highest Found Minimum Distance}
\label{tab:circulant}
\begin{tabular}{cccc}
\toprule
$n$ & Max $d$ & \# Codes & \# Self-dual \\
\midrule
  2 &  2 &    1 &  1 \\
  3 &  2 &    2 &  1 \\
  4 &  3 &    1 &  0 \\
  5 &  3 &    3 &  1 \\
  6 &  4 &    1 &  1 \\
  7 &  4 &    2 &  0 \\
  8 &  4 &   11 &  1 \\
  9 &  4 &   22 &  2 \\
 10 &  5 &    4 &  0 \\
 11 &  5 &   21 &  0 \\
 12 &  6 &    2 &  1 \\
 13 &  6 &    2 &  0 \\
 14 &  6 &   54 &  3 \\
 15 &  6 &  325 &  3 \\
 16 &  7 &    1 &  0 \\
 17 &  7 &    9 &  1 \\
 18 &  8 &    1 &  1 \\
 19 &  7 & 1366 &  4 \\
 20 &  8 &    4 &  3 \\
 21 &  8 &   42 &  0 \\
 22 &  8 & 1328 & 17 \\
 23 &  8 & 8027 &  2 \\
 24 &  9 &    1 &  0 \\
 25 &  9 &   25 &  0 \\
 26 &  9 & 1877 &  0 \\
\bottomrule
\end{tabular}
\end{table}

\begin{proposition}
A circulant or bordered circulant directed graph code will always be isodual.
\end{proposition}
\begin{proof}
A circulant directed graph code of length $n$ has generator matrix
\[
C = 
\begin{pmatrix}
\omega & a_1    & a_2    & \cdots & a_n    \\
a_n    & \omega & a_1    & \cdots & a_{n-1} \\
a_{n-1} & a_n    & \omega & \cdots & a_{n-2} \\ 
\vdots & \vdots & \vdots & \ddots & \vdots \\
a_1    & a_2    & a_3     & \cdots & \omega
\end{pmatrix},
\]
where $(a_1,a_2,\ldots,a_n)$ is any binary sequence of length $n-1$.
It follows from Prop.~\ref{dualtranspose} that the dual code is generated by
\[
C^T = 
\begin{pmatrix}
\omega & a_n    & a_{n-1} & \cdots & a_1     \\
a_1    & \omega & a_n     & \cdots & a_2    \\
a_2    & a_1    & \omega  & \cdots & a_3    \\ 
\vdots & \vdots &  \vdots & \ddots & \vdots \\
a_n    & a_{n-1} & a_{n-2} & \cdots   & \omega
\end{pmatrix}.
\]
We can obtain $C^T$ from $C$ by reversing the order of the columns, and then
reversing the order of the rows. Permuting rows has no effect on the code, and permuting
columns produces an equivalent code. Hence the code must be equivalent to its dual.
The same argument holds for bordered circulant codes, except that the first row and 
column remain fixed.\qed
\end{proof}

If $\mathcal{C}$ is a circulant directed graph code of length $n$, then
$\left|\Aut(\mathcal{C})\right|$ must be divisible by $n$, since the structure of the generator matrix ensures
that cyclically shifting all codewords will preserve the code. Similarly, 
if $\mathcal{C}$ is a bordered circulant directed graph code of length $n$, then
$\left|\Aut(\mathcal{C})\right|$ must be divisible by $n-1$, since the code is preserved
by fixing the first coordinate and cyclically shifting the last $n-1$ coordinates.

With our method, we are able to find new codes, since the existence of near-extremal
formally self-dual codes of lengths 11 and 13 was previously an open problem~\cite{hankim}.
We also answer the question of the existence of an additive $(13,2^{13},6)$ code in the positive.
Parameters for optimal additive codes over $\GF(4)$ were determined for $n \le 12$ by 
Blockhuis and Brouwer~\cite{blokhuis}, and for $n \le 13$ by Bierbrauer~et~al.~\cite{bierbrauer},
by using geometric descriptions of codes.
Bierbrauer~et~al. found an additive $(11,2^{11},5)$ code, but posed as an open question the existence 
of an additive $(13,2^{13},6)$ code.

There are at least 21 formally self-dual $(11,2^{11},5)$ codes, 
with generators available at \url{http://www.ii.uib.no/~larsed/directed/}.
We find codes with five different weight enumerators:
\begin{eqnarray*}
W_{11,1}(1,y) &=& 1 + 55y^5 + 242y^6 + 275y^7 + 495y^8 + 605y^9 + 286y^{10} + 89y^{11}, \\
W_{11,2}(1,y) &=& 1 + 66y^5 + 198y^6 + 330y^7 + 495y^8 + 550y^9 + 330y^{10} + 78y^{11}, \\
W_{11,3}(1,y) &=& 1 + 70y^5 + 182y^6 + 350y^7 + 495y^8 + 530y^9 + 346y^{10} + 74y^{11}, \\
W_{11,4}(1,y) &=& 1 + 75y^5 + 162y^6 + 375y^7 + 495y^8 + 505y^9 + 366y^{10} + 69y^{11}, \\
W_{11,5}(1,y) &=& 1 + 77y^5 + 154y^6 + 385y^7 + 495y^8 + 495y^9 + 374y^{10} + 67y^{11}.
\end{eqnarray*}
There are five codes with weight enumerator $W_{11,1}$. Of these, three have automorphism
groups of order 11, one has 10 automorphisms, and one has 110 automorphisms.
There are four codes with $W_{11,2}$, all with 11 automorphisms.
There are two codes with $W_{11,3}$, both with 10 automorphisms.
There are five codes with $W_{11,4}$, all with 10 automorphisms.
There are five codes with $W_{11,5}$, all with 11 automorphisms.

We have found two formally self-dual $(13,2^{13},6)$ codes.
$\mathcal{C}_{13,1}$ is generated by all cyclic shifts of $(\omega 101001110000)$ and
has an automorphism group of order $13$.
$\mathcal{C}_{13,2}$ is generated by all cyclic shifts of $(\omega 111011111010)$ and
has an automorphism group of order $78$.
Both these codes have the same weight enumerator:
\[
\begin{split}
W_{13}(1,y) &= 1 + 247y^6 + 481y^7 + 936y^8 + 1625y^9 + 2197y^{10} + 1755y^{11} + \\
  &\quad 715y^{12} + 235y^{13}.
\end{split}
\]

Note that for several lengths, there are no self-dual codes among the circulant and bordered circulant 
codes with highest minimum distance. The best known self-dual codes of length 24, 25, and 26 have minimum distance 8.
We find a single isodual $(24,2^{24},9)$ code generated by the cyclic shifts of 
$(\omega 01101111111111010000110)$ with automorphism group of order $72$ and weight enumerator
\[
\begin{split}
W_{24}(1,y) &= 1 + 1752y^9 + 8748y^{10} + 26064y^{11} + 81408y^{12} + 232776y^{13} + \\
  &\quad 573516y^{14} + 1119264y^{15} + 1869777y^{16} + 2676456y^{17} + 3096804y^{18} + \\
  &\quad 2959056y^{19} + 2204568y^{20} + 1255416y^{21} + 520740y^{22} + 134208y^{23} + \\
  &\quad 16662y^{24}.
\end{split}
\]
We also find 25 isodual $(25,2^{25},9)$ codes with 25 different weight enumerators, 
and 1877 isodual $(26,2^{26},9)$ codes with 1865 different weight enumerators.

We have previously studied circulant undirected graph codes~\cite{weighted}. There are only 
$2^{\left\lceil\frac{n-1}{2}\right\rceil}$ such codes of length $n$, due to the fact that the
generator matrix must be symmetric.
Gulliver and Kim~\cite{gulliverkim} also performed a computer search of circulant self-dual
additive codes over $\GF(4)$, but their search was not restricted to 
graph codes. 

A particularly interesting type of circulant code is a type of \emph{quadratic residue code}~\cite{handbook}.
The length of such a code must be a prime $p$. When $p \equiv 1 \pmod{4}$, the quadratic residue 
code will be self-dual, and the corresponding undirected graph is known as a \emph{Paley graph}.
When $p \equiv 3 \pmod{4}$, the code will only be isodual.
The first row of the generator matrix of the code is $(\omega, l_1, \ldots, l_{p-1})$, where
$l_i = 1$ if $i$ is a quadratic residue modulo $p$, i.e., if $x^2 \equiv i \pmod{p}$ has a solution
$x \in GF(p)$. Otherwise, $l_i = 0$. Many codes with high minimum distance
can be obtained from this construction. For instance, by bordering quadratic residue codes,
as described above, we obtain self-dual $(6,2^6,4)$, $(14,2^{14},6)$, and $(30,2^{30},12)$ codes,
and isodual $(4,2^4,3)$, $(8,2^8,4)$, and $(12,2^{12},6)$ codes.

\paragraph*{Acknowledgement}
The authors would like to thank Jürgen Bierbrauer for helpful comments.
This research was supported by the Research Council of Norway.

\end{document}